\documentclass[10pt]{article}
\usepackage{amssymb}
\usepackage{amsmath}
\usepackage{amsthm}

\newtheorem{theorem}{Theorem}

\newtheorem{proposition}{Proposition}

\newtheorem{lemma}{Lemma}

\title
{ Determination of a source term in Schr\"odinger equation with a data taken at final moment of observation.}

\author{
O.Y. Imanuvilov and $^{2,3}$ M.~Yamamoto
}

\date{}
\begin{document}
\maketitle
\begin{abstract}{\it We establish the Lipshitz stability estimate in inverse problem of determination of a source term or zero order  term in the Schr\"odinger equation  with time-dependent coefficients under some non-trapping assumption. Based on this result we established the Lipshitz stability of the determination of  a real-valued coefficient corresponding to zero-th order term in the Schr\"odinger equation.}
\end{abstract}


\section{Introduction}
The paper is concerned  with two inverse problems for the Sch\"odinger equation: the determination of  a source term in the right hand side of the Sch\"odinger equation and determination of  coefficient respect to the zero order of the  Sch\"odinger equation. The data is collected at the final moment of observation and on a part of a boundary.
In the situation when a part of a boundary satisfies some geometric condition (more specifically  non-trapping condition for   bicharacteristcs associated with elliptic part  of Schr\"odinger operator) we prove the  Lipshitz stability for determination of the source term and the Lipshitz stability for determination of the zero order coefficient in the Schr\"odinger equation.

The uniqueness in determining the $L^\infty$-potential in the  Schr\"odinger equation was firstly
proved by Bukhgeim \cite{Buh} for the case of one dimensional Schr\"odinger equation.  After that, the uniqueness and Lipshitz stability of   determination of the source term in the right hand side of the  Schr\"odinger equation was studied in \cite{B-P}, \cite{BP1}  and \cite{MOR} for the case when  the principal term in the elliptic part of the Schr\"odinger equation is the Laplace operator. In \cite{B-P} and \cite{MOR} the Lipshitz stability for determination of the source  term was established under condition that trace of solution at final moment of observation is either  real valued function or purely imaginary function. One of goals of this work is to remove this assumption.
It should be mentioned that in  \cite {MOR} the main machinery used by the authors are Carleman estimates.
In \cite{BM} the authors studied an inverse problem of determination of coefficient of zero order term for the Schr\"odinger equation with piecewise constant coefficient in principal part.
In 
\cite {KYM}  authors proves Lipshitz stability in the determination of magnetic potential vector for the magnetic Schr\"odinger equation.

In 
\cite {YM} the authors established  observability estimates for the Schr\"odinger equation with integral terms.

The proof of the conditional Lipshitz stability stability is based on the  Carleman estimates  with boundary. In that sense we are following  the pioneering work of M. Klibanov and A. Bukhgeim \cite{BK}. 
  
Unlike \cite{B-P} we do not use an extension of the solution to the Schr\"odinger equation. As a result of this  our weight function is  not degenerates near the final moment of  observation. This results in additional terms  in the  right hand side of Carleman estimate, containing the trace of function at final moment of observation. 
Finally we mention that uniqueness  of determination of the source term in the  right hand side of the Schr\"odinger equation without geometrical restrictions on  a domain of observation was established in  \cite {IY3} for the right hand side of the form $R(t) f(x).$


\section{Inverse Problem}
Let $\Omega$ be a bounded domain in $\Bbb R^n$ with $C^2-$ boundary. Let $\Gamma$ be an arbitrary fixed subboundary of $\partial\Omega.$We set 
$Q=(0,T)\times \Omega,\Sigma=(0,T)\times \partial\Omega,\Sigma_1=(0,T)\times \partial\Omega\setminus \bar\Gamma,\Sigma_0=(0,T)\times \Gamma.$
Let function $u$ solve the following boundary value problem
\begin{equation}\label{Ssok1}
P(t,x,D)u=i\partial_tu-\sum_{{\ell},j=1}^n \partial_{x_{\ell}}(a_{{\ell}j}(t,x)\partial_{x_j} u)+\sum_{{\ell}=1}^n b_{\ell}(t,x)\partial_{x_{\ell}} u+c(t,x) u=R(t,x) f(x)\quad \mbox{in}\quad Q,
\end{equation}
\begin{equation}\label{Ssok2}
u\vert_{\Sigma}=0,
\end{equation}\begin{equation}\label{Ssok3}
u(T,\cdot)=u_T.
\end{equation}
Complex-valued functions $u_T$ and $R$ are given. A real valued source term $f$ is unknown.

 We consider the following inverse source problem:

 {\it Suppose that $\partial_{\nu} u\vert_{\Sigma_0}$ is given. Determine the source term $f.$}

We will make the following assumptions: functions $a_{{\ell}j}$ are real-valued  for all ${\ell},j$ from $\{1,\dots,n\}$ and 
\begin{equation}\label{sok4}
a_{{\ell}j}=a_{j\ell} \in C^{1}(\bar Q), \quad b_{\ell}\in C^{1,0}(\bar Q),\quad c\in  C^{1,0}(\bar Q)\quad \forall {\ell},j\in \{1,\dots,n\}
\end{equation}
and  there exist a constant $\beta>0$ such that
\begin{equation}\label{sok5}
\sum_{{\ell},j=1}^n a_{{\ell}j}(t,x)\eta_{\ell}\eta_j\ge \beta\vert \eta\vert^2\quad \forall (t, x,\eta)\in Q\times \Bbb R^n.
\end{equation}
In addition assume that there exists  strictly positive constant $\beta_1$ such that
\begin{equation}\label{sok6}
R\in C^{2}(\bar Q), \quad \vert R(T,x)\vert >\beta_1>0\quad \mbox{on }\,\, \bar \Omega.
\end{equation}

Denote  by $p(t,x,\xi)$ the principal symbol of the Schr\"odinger operator $P:$
$$
p(t,x,\xi)=-\xi_0+a(t,x,\xi',\xi'),\quad a(t,x,\xi',\eta)=\sum_{k,j=1}a_{kj}(t,x)\xi_k\eta_j,\quad 
$$
$$\xi=(\xi_0,\xi'), \,\, \xi'=(\xi_1,\dots ,\xi_n), \eta=(\eta_1,\dots,\eta_n).
$$
We set $M_\lambda \xi=(\lambda^2\xi_0, \lambda\xi_1,\dots \lambda\xi_n).$
For any function $\tilde \psi\in C^{1,2}(\bar Q)$ we introduce the function
$q_{\tilde \psi}(t,x,\xi,\tau)$:
on the set $\{(t,x,\xi,\tau): (t,x)\in Q,\quad \vert\xi_0\vert+\vert \xi'\vert^2+\tau^2=1\}$
$$
q_{\tilde \psi}(t,x,\xi,\tau)=\lim_{\lambda \rightarrow +\infty}\frac{\{p(t,x,M_\lambda \xi-i\lambda \tau\nabla_x\tilde \psi),p(t,x, M_\lambda \xi+i\lambda \tau\nabla_x \tilde\psi)\}}{2i\lambda \tau },
$$
where
$$
\nabla_x \tilde \psi=(0,\partial_{x_1}\tilde \psi,\dots,\partial_{x_n}\tilde\psi).
$$

On the set $\{(t,x,\xi,\tau): (t,x)\in Q,\quad (\xi,\tau)\ne 0\}$
$$
q_{\tilde \psi}(t,x,\xi,\tau)=q_{\tilde \psi}(t,x,\frac{\xi_0}{{\vert\xi_0\vert+\vert \xi'\vert^2+\tau^2}},\frac{\xi_1}{\root\of{\vert\xi_0\vert+\vert \xi'\vert^2+\tau^2}},\dots,\frac{\xi_n}{\root\of{\vert\xi_0\vert+\vert \xi'\vert^2+\tau^2}},\frac{\tau}{\root\of{\vert\xi_0\vert+\vert \xi'\vert^2+\tau^2}} )
$$

In order to formulate a Carleman estimate we introduce the notion of pseudo-convex function respect to the  principal symbol of the  operator $P:$

{\bf Definition.} {\it 
A function $\tilde \psi\in C^2(\bar Q),\nabla_{x} \tilde \psi\ne 0$ on $\bar Q$  is called pseudo-convex respect to the symbol $p$ if 
\begin{eqnarray}\label{pseudo}
q_{\tilde \psi}(t,x,\xi,\tau)>0\quad\mbox{on}\quad \{ (t,x,\xi,\tau)\in Q\times \Bbb R^{n+1}\times \Bbb R_+ , p(t,x,\xi-i\tau\nabla_x\tilde \psi)=0,\quad \vert \xi_0\vert+\vert \xi'\vert^2+\tau^2=1,\tau>0\}.
\end{eqnarray}}
Later we will introduce the Carleman estimate with boundary. Therefore, except a standard assumption of pseudoconvexity of a weight function, we introduce additional assumption on a weight function on  an unobserved part of the boundary:

{\bf Condition 1.} {\it There exist a function $\psi$ pseudo-convex respect to the principal symbol of the operator $p$ such that
\begin{equation}\label{sokol}
a(t,x,\nu,\nabla \psi)<0\quad \mbox {on}\quad [0,T]\times  (\partial\Omega\setminus\Gamma).
\end{equation}}

We set for any $\mbox{\bf v},\mbox{\bf w}\in \Bbb R^n$ and $p\in \{0,1,\dots,n\}$ we define the quadratic forms $a_{(p)}$ by the formula:
\begin{equation}\label{pioner}
a_{(0)}(t,x,\mbox{\bf v},\mbox{\bf w})=\sum_{{\ell},j=1}^n \partial_{t}a_{{\ell}j}(t,x)v_{\ell}w_j,\quad a_{(p)}(t,x,\mbox{\bf v},\mbox{\bf w})=\sum_{{\ell},j=1}^n \partial_{x_p}a_{{\ell}j}(t,x)v_{\ell}w_j.
\end{equation}
We have
\begin{proposition}\label{11}
Let function $\tilde \psi\in C^2(\bar Q),\nabla_{x} \tilde \psi\ne 0$ on $\bar Q$ satisfies the inequality
\begin{equation}\label{ramsai}4\sum_{m,\tilde m,k,\ell=1}^n  a_{mk}a_{\tilde m \ell}(\xi_m\xi_{\tilde m}) \partial^2_{x_kx_\ell}\tilde \psi
-2\sum_{m,k=1}^na_{mk}\partial_{x_m}\tilde \psi a_{(k)}(t,x,\xi',\xi')
+4\sum_{m,k=1}^na_{mk}\xi_m a_{(k)}(t,x,\xi',\nabla\tilde \psi))>0\end{equation}
for all $(t,x,\xi')\in \{(t,x,\xi')\in Q\times \Bbb R^n\vert a(t,x,\nabla\tilde \psi,\xi')=0\}.$
Then for all sufficiently large $\lambda$ function $e^{\lambda\tilde \psi}$ is pseudo-convex respect to the principal symbol of operator $p$.
\end{proposition}
Proof of Proposition \ref{11} is given in section 3.

{\bf Remark.}
{\it One can check by direct computations that if function $\psi$ is pseudo-convex the inequality (\ref{ramsai}) holds true. If $a_{\ell j}=\delta_{\ell j}$ function $\psi(x)=\vert x-y\vert^2, y\notin \bar\Omega$ satisfies (\ref{ramsai}).}

Our main result is the following
\begin{theorem}\label{SOKOL} Let (\ref{sok4})-(\ref{sok6}) and Condition 1 holds true. Let  $u,\partial_t u\in H^{1,2}(Q), u(T,\cdot)\in H^3(\Omega)\cap H^1_0(\Omega),  f\in L^2(\Omega)$  and $\tilde u,\partial_t \tilde u\in H^{1,2}(Q), \tilde u(T,\cdot)\in H^3(\Omega)\cap H^1_0(\Omega), \tilde f\in L^2(\Omega)$ are both solutions to problem (\ref{Ssok1})-(\ref{Ssok3}) and functions $f,\tilde f$ are real-valued.  Then there exist a constant $C$ independent of $u,f,\tilde u,\tilde f$ such that
\begin{equation}\label{pobeda}
\Vert f-\tilde f\Vert_{L^2(\Omega)}\le C(\Vert u(T,\cdot)-\tilde u(T,\cdot)\Vert_{H^3(\Omega)}+\sum_{k=0}^1\Vert\partial^k_t \partial_\nu( u-\tilde u)\Vert_{L^2(\Sigma_0)}).
\end{equation}
\end{theorem}

{\bf Proof.}  Proof of the Theorem \ref{SOKOL} is based on a Carleman estimate. First we introduce weight functions which will be used  to set up such a  Carleman estimate.

Using function $\psi$ from Condition 1 we introduce two functions $\alpha(t,x)$ and $\varphi(t,x)$:
\begin{equation}\label{DD1.2} \varphi(t, x) = e^{\lambda \psi} / t,
\quad 
 \alpha(t, x) = (e^{\lambda\psi}
- e^{2\lambda ||\psi||_{C^0(\overline{\Omega})}}) / t.
\end{equation}

Let us consider the boundary value problem
\begin{equation}\label{LD1.4}
 P (t,x,D)z =  g \quad \mbox{in}\quad Q,
\end{equation}
\begin{equation}\label{LD1.5}z \big|_\Sigma = 0,\quad z(0, \cdot) = z_0.
\end{equation}

Later we will show that   function $\alpha(t,x)$ is also pseudo-convex respect to  the symbol $p.$

In order to formulate the next lemma  we introduce two operators

\begin{equation}\label{sosna}
P_1(t,x,D,\tau)u=2\tau\varphi\lambda a(t,x,\nabla \psi,\nabla u)+\tau\lambda^2\varphi a(t,x,\nabla\psi,\nabla\psi)u+\tau\varphi\lambda \sum_{k,j=1
}^n \partial_{x_k} a_{kj}\partial_{x_j}\psi u ,
\end{equation}
$$
P_2(t,x,D,\tau)u=i\partial_t u- \sum_{{\ell},j=1}^n a_{{\ell}j} \partial^2_{x_{\ell}x_j} u -
\lambda^2 {\tau}^2 \varphi^2 a(t, x, \nabla\psi, \nabla\psi) u.
$$
Let $\gamma=\sum_{j=1}^n\Vert b_j\Vert_{L^\infty(Q)}+\Vert c\Vert_{L^\infty(Q)}.$ In the space $L^2(G)$ we introduce the scalar product by formula$$
(u,v)_{L^2(G)}=\int_G u\bar v dy.$$
We have the following:

\begin{lemma}\label{D1.2} {\it Let 
$$
a_{{\ell}j}=a_{j\ell} \in C^{1}(\bar Q), \quad b_{\ell}\in L^\infty(\bar Q),\quad c\in  L^\infty(\bar Q)\quad \forall {\ell},j\in \{1,\dots,n\},
$$
 (\ref{sok5}) is fulfilled and functions $\varphi, \alpha$
 be defined as in (\ref{DD1.2}), where function $\psi$ satisfies Condition 1. Let $z\in H^{1,2}(Q)$ and $g\in L^2(Q).$ Then there exists  constants
$\hat\lambda(\gamma) > 0$ such that for an arbitrary $\lambda \ge \hat{\lambda}$ there exists
  ${\tau}_0(\lambda,\gamma)$ such that for each ${\tau} \ge \tau_0(\lambda,\gamma)$ the solutions of problem (\ref{LD1.4}) - (\ref{LD1.5})
satisfy the following inequality:
\begin{eqnarray}\label{DD1.6}
  \int_Q ( {\tau}\varphi |\nabla z|^2 + {\tau}^3 \varphi^3 \vert z\vert^2)
e^{2{\tau}\alpha}   dx\,dt+\sum_{k=1}^2\Vert P_k(t,x,D,\tau) (ze^{\tau\alpha})\Vert^2_{L^2(Q)}\nonumber\\
 \le C_1\big(\vert \mbox{Re}\, (i (ze^{\tau\alpha})(T,\cdot),P_1(T,x,D,\tau)(ze^{\tau\alpha}) (T,\cdot))_{L^2(\Omega)}\vert+ \int_Q |g|^2 e^{2{\tau}\alpha}  dxdt +
\int_{\Sigma_0} \tau \varphi \vert\partial_\nu z\vert^2 
e^{2{\tau}\alpha}  d\Sigma\big).
\end{eqnarray}
}
\end{lemma}
Now we proceed to the proof of the theorem \ref{SOKOL}. without loss of generality by (\ref{sok6}) we may assume that 
\begin{equation}\label{Asok6}
R\in C^{2}(\bar Q), \quad \vert R(t,x)\vert >\beta_2>0\quad \mbox{on }\,\, \bar Q.
\end{equation}
Let pairs $(u,f)$ and $(\tilde u,\tilde f)$ be solutions to problem (\ref{Ssok1})-(\ref{Ssok3}). The pair $(w,q)=(u-\tilde u,f-\tilde f)$ satisfies 
\begin{equation}\label{mob}
P(t,x,D) w=Rq\quad\mbox{in}\,\, Q, \quad w(T,\cdot)=u(T,\cdot)-\tilde u(T,\cdot),\quad w\vert_\Sigma=0,\quad  \partial_\nu w\vert_{\Sigma_0}=\partial_\nu(u-\tilde u).\end{equation}

By (\ref{mob}) the 
function $\tilde w=i\bar R w$ solves the initial value problem
\begin{equation}\label{slon}
\tilde P(t,x,D)\tilde w=i\vert R\vert^2 q \quad\mbox{in}\,\, Q, \quad \tilde w(T,\cdot)=i\bar R(u(T,\cdot)-\tilde u(T,\cdot)),\quad \tilde w\vert_\Sigma=0,\quad  \partial_\nu \tilde w\vert_{\Sigma_0}=i\bar R\partial_\nu(u-\tilde u).
\end{equation}
Operator $\tilde P$ is given by formula
$$
\tilde P(t,x,D)r =P(t,x,D)r+\frac{1}{\bar R}a(t,x,\nabla \bar R,\nabla r)-\frac{1}{\bar R}a(t,x,\nabla \bar R,\nabla \bar R)r-\sum_{j=1}^n b_j\partial_{x_j}\bar R r/\bar R+\sum_{{\ell},j=1}^n a_{{\ell}j}\partial^2_{x_{\ell}x_j}\bar R r/\bar R=
$$
$$
i\partial_tr-\sum_{{\ell},j=1}^n \partial_{x_{\ell}}(a_{{\ell}j}(t,x)\partial_{x_j} r)+\sum_{{\ell}=1}^n\tilde  b_{\ell}(t,x)\partial_{x_{\ell}} r+\tilde c(t,x) r.
$$
The principal symbols of operators $P$ and $\tilde P$ are the same. So Condition 1 holds true with the same function $\psi$ for the operator $\tilde P.$  By (\ref{sok4}), (\ref{Asok6}) the coefficients respect to the first and zero order terms of the operator $\tilde P$ satisfy (\ref{sok4}).

Hence, 
differentiating equation (\ref{slon}) respect to $t$ we obtain

   \begin{equation}\label{silos}
\tilde P(t,x,D)\partial_t\tilde w=i\partial_t\vert R\vert^2 q+\mathcal M(t,x,D)\tilde w\quad\mbox{in}\,\, Q, \quad \partial_t\tilde w(T,\cdot)=-\vert R(T,\cdot)\vert^2 q +\mathcal K(u(T,\cdot)-\tilde u(T,\cdot)),\quad \partial_t\tilde w\vert_\Sigma=0,
\end{equation} where 
$$
\mathcal K \tilde w=i(\sum_{{\ell},j=1}^n \partial_{x_{\ell}}(a_{{\ell}j}(T,x)\partial_{x_j} \tilde w)-\sum_{{\ell}=1}^n\tilde  b_{\ell}(T,x)\partial_{x_{\ell}} \tilde w-\tilde c(T,x) \tilde w),
$$
$$
\mathcal M(t,x,D)\tilde w=\sum_{{\ell},j=1}^n \partial_{x_{\ell}}(\partial_ta_{{\ell}j}(t,x)\partial_{x_j} \tilde w)-\sum_{{\ell}=1}^n\partial_t\tilde  b_{\ell}(t,x)\partial_{x_{\ell}} \tilde w-\partial_t\tilde c(t,x) \tilde w.
$$
Applying to equation (\ref{silos}) Carleman estimate (\ref{DD1.6}) we obtain that there exist a constant $C_1$ and $\tau_0$ such that
\begin{eqnarray}\label{silos1}
  \int_Q ( {\tau}\varphi |\nabla\partial_t \tilde w|^2 + {\tau}^3 \varphi^3\vert\partial_t \tilde w\vert^2)
e^{2{\tau}\alpha}   dx\,dt+\sum_{k=1}^2\Vert P_k(t,x,D,\tau) (\partial_t\tilde we^{\tau\alpha})\Vert^2_{L^2(Q)}\nonumber\\
 \le C_1\big(\vert \mbox{Re}\, (i (\partial_t\tilde w e^{\tau\alpha})(T,\cdot),P_1(T,x,D,\tau)(\partial_t\tilde we^{\tau\alpha}) (T,\cdot))_{L^2(\Omega)}\vert\nonumber\\+ \int_Q (|q|^2+\vert \mathcal M\tilde w\vert^2) e^{2{\tau}\alpha}  dxdt +
\int_{\Sigma_0} \tau \varphi \vert\partial_\nu\partial_t (u-\tilde u)\vert^2 
e^{2{\tau}\alpha}  d\Sigma\big)
\end{eqnarray}  for all $\tau\ge\tau_0.$
Observe that by (\ref{silos}) 
$$
 \mbox{Im}\, ( (\partial_t\tilde w e^{\tau\alpha})(T,\cdot),P_1(T,x,D,\tau)(\partial_t\tilde we^{\tau\alpha}) (T,\cdot))_{L^2(\Omega)}=
 \mbox{Im}\, ( (\mathcal K \tilde we^{\tau\alpha})(T,\cdot),P_1(T,x,D,\tau)(\mathcal K \tilde we^{\tau\alpha}) (T,\cdot))_{L^2(\Omega)}
$$
$$+\mbox{Im}\, ( -\vert R(T,\cdot)\vert^2 q,P_1(T,x,D,\tau)(\mathcal K \tilde we^{\tau\alpha}) (T,\cdot))_{L^2(\Omega)}
$$
$$
-\mbox{Im}\, ( (\mathcal K \tilde we^{\tau\alpha})(T,\cdot),P_1(T,x,D,\tau)\vert R(T,\cdot)\vert^2 q)_{L^2(\Omega)}.
$$
From this equality we obtain
\begin{eqnarray}\label{ox}
\vert \mbox{Re}\, (i (\partial_t\tilde w e^{\tau\alpha})(T,\cdot),P_1(T,x,D,\tau)(\partial_t\tilde we^{\tau\alpha}) (T,\cdot))_{L^2(\Omega)}\vert\nonumber\\\le C_2(\tau\Vert  \tilde w(T,\cdot)e^{\tau\alpha(T,\cdot)}\Vert^2_{H^{3,\tau}(\Omega)}+\frac 1\tau\Vert qe^{\tau\alpha(T,\cdot)}\Vert_{L^2(\Omega)}^2).
\end{eqnarray}
Here 
$$\Vert v\Vert_{H^{3,\tau}(\Omega)}=\Vert v\Vert_{H^{3}(\Omega)}+\tau^3\Vert v\Vert_{L^2(\Omega)}.$$
Using the inequality (\ref{ox}) in the estimate (\ref{silos1}) we have
\begin{eqnarray}\label{prison2}
 \int_Q ( {\tau}\varphi |\nabla\partial_t \tilde w|^2 + {\tau}^3 \varphi^3\vert \partial_t\tilde w\vert^2)
e^{2{\tau}\alpha}   dx\,dt+\sum_{k=1}^2\Vert P_k(t,x,D,\tau) (\partial_t\tilde we^{\tau\alpha})\Vert^2_{L^2(Q)}\le C_3\big( \frac 1\tau\Vert qe^{\tau\alpha(T,\cdot)}\Vert_{L^2(\Omega)}^2
\nonumber\\
 +\int_Q \vert \mathcal M\tilde w\vert^2 e^{2{\tau}\alpha}  dxdt +\tau\Vert  \tilde w(T,\cdot)e^{\tau\alpha(T,\cdot)}\Vert^2_{H^{3,\tau}(\Omega)}+
\int_{\Sigma_0} \tau \varphi \vert\partial_\nu \partial_t\tilde w\vert^2 
e^{2{\tau}\alpha}  d\Sigma\big)\quad \forall\tau\ge\tau_1.
\end{eqnarray}
Now we need to estimate $\int_Q \vert \mathcal M\tilde w\vert^2 e^{2{\tau}\alpha}  dxdt $. In order to do that we write equation (\ref{slon}) as
$$
-\sum_{{\ell},j=1}^n \partial_{x_{\ell}}(a_{{\ell}j}(t,x)\partial_{x_j} \tilde w)=-\sum_{{\ell}=1}^n b_{\ell}(t,x)\partial_{x_{\ell}} \tilde w-c(t,x) \tilde w-i\partial_t\tilde w+i\vert R\vert^2 q \quad \mbox{in}\quad Q.
$$
For each $t\in [0,T]$ we apply to this equation a Carleman estimate  for elliptic operator (see e.g. \cite{Im})
\begin{equation}\label{prison}
\int_\Omega(\frac{1}{\tau \varphi} \sum_{\ell,k=1}^n\vert \partial^2_{x_\ell x_j} \tilde w\vert^2+\tau\varphi\vert \nabla w\vert^2+\tau^3\varphi^3\vert \tilde w\vert^2)e^{2\tau\alpha}dx \le C_4(\int_\Omega (\vert q\vert^2 +\vert\partial_t\tilde w\vert^2)e^{2\tau\alpha}dx+\int_{\Gamma}\tau\varphi \vert\partial_\nu \tilde w\vert^2 d\sigma)\quad \forall \tau\ge\tau_2.
\end{equation}
Integrating the inequality (\ref{prison}) on time interval $(0,T)$ we have
\begin{eqnarray}\label{prison1}
\int_Q(\sum_{\ell,k=1}^n\vert \partial^2_{x_\ell x_j} \tilde w\vert^2+\tau^2\varphi^2\vert \nabla w\vert^2+\tau^4\varphi^4\vert \tilde w\vert^2)e^{2\tau\alpha}dx \le C_5(\int_\Omega \vert q\vert^2e^{2\tau\alpha(T,x)}dx\nonumber\\+\int_Q\tau\varphi\vert\partial_t\tilde w\vert^2e^{2\tau\alpha}dxdt+\int_{\Sigma_0}\tau^2\varphi^2 \vert\partial_\nu \tilde w\vert^2 d\Sigma) \quad \forall \tau\ge\tau_2.
\end{eqnarray}
Hence  by (\ref{prison1}) and (\ref{prison2}) we have
\begin{eqnarray}\label{silon4}
\Vert P_2 (\partial_t\tilde w e^{s\alpha})\Vert_{L^2(Q)}^2 +\tau^3\Vert \partial_t\tilde we^{\tau\alpha}\Vert^2_{L^2(Q)}\le  C_6(\Vert qe^{\tau\alpha(T,\cdot)}\Vert^2_{L^2(\Omega)} \nonumber\\
+\Vert  \tilde w(T,\cdot)e^{\tau\alpha(T,\cdot)}\Vert^2_{H^{3,\tau}(\Omega)}+\int_{\Sigma_0} \tau \varphi \vert\partial_\nu\partial_t \tilde w\vert^2 
e^{2{\tau}\alpha}  d\Sigma+\int_{\Sigma_0} \tau \varphi \vert\partial_\nu\partial_t \tilde w\vert^2 
e^{2{\tau}\alpha}  d\Sigma)\quad \forall \tau\ge\tau_3.
\end{eqnarray}
Since 
$$
\tau^\frac 34 \Vert \partial_t\tilde w(T,\cdot)e^{\tau\alpha(T,\cdot)}\Vert^2_{L^2(\Omega)}\le C_3(\Vert P_2(\partial_t\tilde w e^{s\alpha})\Vert^2_{L^2(Q)}+\tau^3\Vert \partial_t\tilde w e^{s\alpha}\Vert^2_{L^2(Q)}) \quad \forall \tau\ge\tau_4 .
$$
Form (\ref{silon4}) and (\ref{silos}) we have
$$
\tau^\frac 34 \Vert \partial_t\tilde w(T,\cdot)e^{\tau\alpha(T,\cdot)}\Vert^2_{L^2(\Omega)}=\tau^\frac 34 \Vert (\vert R\vert^2  (T,\cdot) q-\mathcal K\tilde w(T,\cdot))e^{\tau\alpha(T,\cdot)}\Vert^2_{L^2(\Omega)}
$$
$$\le C_7(\Vert P_2(\partial_t\tilde w e^{s\alpha})\Vert^2_{L^2(Q)}+\tau^\frac 32\Vert \partial_t\tilde w e^{\tau\alpha}\Vert^2_{L^2(Q)})
$$
$$\le C_8(\Vert qe^{\tau\alpha(T,\cdot)}\Vert^2_{L^2(\Omega)}+\Vert  \tilde w(T,\cdot)e^{\tau\alpha(T,\cdot)}\Vert^2_{H^{3,\tau}(\Omega)}+\int_{\Sigma_0} \tau \varphi \vert\partial_\nu\partial_t \tilde w\vert^2 
e^{2{\tau}\alpha}  d\Sigma+\int_{\Sigma_0} \tau \varphi \vert\partial_\nu\partial_t \tilde w\vert^2 
e^{2{\tau}\alpha}  d\Sigma).
$$
for all $\tau\ge\tau_5.$ 
Taking  parameter $\tau$ sufficiently large we obtain (\ref{pobeda}). Proof of the theorem is complete.
$\blacksquare$

Next we consider the inverse problem of determination of the coefficient of zero order term in the parabolic equation in the case  when data is given in the initial time moment
\begin{equation}\label{PZSsok1}
P_j(t,x,D)u_j=i\partial_tu_j-\sum_{\ell,k=1}^n \partial_{x_\ell}(a_{\ell k}(t,x)\partial_{x_k} u_j)+\sum_{k=1}^n b_k(t,x)\partial_{x_k} u_j+c_j(x) u_j=0\quad \mbox{in}\quad Q,
\end{equation}
\begin{equation}\label{PZSsok2}
 u_j\vert_{\Sigma}=v,
\end{equation}\begin{equation}\label{PZSsok3}
u_j(T,\cdot)=u_{T,j},
\end{equation}
$j\in \{1,2\}.$
Functions $u_{T,j}$ are given. Consider the following inverse problem: {\it Let 
\begin{equation}\label{PP}
\partial_\nu u_1=\partial_\nu u_2\quad\mbox{on} \quad (0,T)\times\Gamma.
\end{equation}
Determine  coefficients  unknown real-values coefficients $c_j.$}

We have
\begin{theorem} Let (\ref{sok4})-(\ref{sok6}) and Condition 1 holds true. Let $c_1(x)$ and $c_2(x)$ be real-valued coefficients.  Assume that $u_T\in H^3(\Omega)$ and $u,\partial_t u\in H^{1,2}(Q).$
Moreover assume that for some positive constants $\beta$ and $M_1, M_2$ 
\begin{equation}\label{meteorit}
\vert u_{T,j}(x)\vert>\beta>0\quad \mbox{on}\,\,\bar \Omega, \quad \Vert u_j\Vert_{C^2(\bar Q)}\le M_1,\quad \Vert c_j\Vert_{L^\infty(\Omega)}\le M_2\quad \forall j\in\{1,2\}.
\end{equation} Then there exist a constant $C$ independent of $u_1,u_2$ such that
\begin{equation}\label{hook}
\Vert c_1-c_2\Vert_{L^2(\Omega)}\le C(\Vert u_{T,1}-u_{T,2}\Vert_{H^3(\Omega)}+\sum_{k=0}^1\Vert\partial^k_t( \partial_\nu u_1-\partial_\nu u_2)\Vert_{L^2(\Sigma_0)}).
\end{equation}
\end{theorem}
{\bf Proof.}
Setting $w=u_1
-u_2, $ we have 
\begin{equation}
P_1(t,x,D)w=u_2(t,x) (c_1-c_2)\quad \mbox{in}\,\, (0,T)\times \Omega, \quad w\vert_{(0,T)\times \partial\Omega}=0, \quad  \quad w(T,\cdot)=u_{T,1}-u_{T,2}.
\end{equation}
There exist a  constant $\kappa$ such that
$\mbox{max}\{\gamma(b_1,\dots,b_n,c_1), \gamma(\partial_tb_1,\dots,\partial_t b_n,c_1)\}\le \kappa$ for all $b_j,c_1$ satisfying (\ref{sok4}) and (\ref{meteorit}).
By (\ref{meteorit}), without loss of generality we may assume that there exist  $\tilde \beta>0$ such that 
\begin{equation} \vert u_2(t,x)\vert>\tilde\beta>0\quad \forall (t,x)\in Q,\quad u_2\in C^2(\bar Q).\end{equation}

Finally we apply Theorem \ref{SOKOL} with $u=w, f=c_2-c_1, R=u_2,\tilde u=0,\tilde f=0$ to obtain estimate (\ref{hook}).
 Proof of the theorem is complete.
$\blacksquare$

 \section{Carleman estimate for the Schr\"odinger equation}

Before starting a proof of  Carleman estimate we establish the following technical proposition.
We set 
$$
P_{1,(0)}(t,x,D,\tau)u=2\tau\varphi\lambda a_{(0)}(t,x,\nabla \psi,\nabla u)+\tau\lambda^2\varphi a_{(0)}(t,x,\nabla\psi,\nabla\psi)u+\tau\varphi\lambda \sum_{k,j=1
}^n \partial_{x_k} \partial_ta_{kj}\partial_{x_j}\psi u ,
$$

\begin{proposition}\label{osa} For any function $w\in H^{1,2}(Q)$ the following equality is true
\begin{eqnarray}\label{suka1}
\mbox{Re} (P_1w,P_2w)_{L^2(Q)}=
\mbox{Re}\, ( w,P_{1,(0)}(t,x,D,\tau)w)_{L^2(Q)}+
\mbox{Re}\, (i w(T,\cdot),P_1(t,x,D,\tau)w (T,\cdot))_{L^2(\Omega)} \nonumber\\
+ \mbox{Re}\, (w, 2{\tau}\lambda i\varphi a(t,x, \nabla\psi, \nu) \partial_tw)_{L^2(\Sigma)}
\nonumber\\
   +\mbox{Re}\int_Q \left(\sum_{{\ell},j=1}^n \partial_{x_j} a_{{\ell}j}
\partial_{x_{\ell}} w P_1 (t,x,D,\tau)\bar w\right. + 2{\tau}\lambda^2 \varphi
 \vert a(t, x, \nabla\psi, \nabla w\big)\vert^2 \nonumber\\
\left.
  +2{\tau}\lambda \varphi \sum_{{\ell},j=1}^n a_{{\ell}j} \partial_{x_{\ell}} w
\sum_{k,\ell=1}^n \partial_{x_j}  (a_{k\ell}\partial_{x_k} \psi)
\partial_{x_\ell} \bar w - {\tau}\lambda\varphi \sum_{k,\ell=1}^n
a_{k\ell} \partial_{x_k}\psi \sum_{{\ell},j=1}^n \partial_{x_\ell} a_{{\ell}j}
\partial_{x_{\ell}} w \partial_{x_j} \bar w \right) dxdt\nonumber \\
 +\mbox{Re} \int_\Sigma \left(
2{\tau}\lambda\varphi 
a(t,x,\nu,\nabla w)a(t,x,\nabla\psi,\nabla\bar w) -
{\tau}\lambda\varphi          a(t, x, \nabla w, \nabla\bar w)
a(t, x, \nu, \nabla\psi)\right) d\Sigma \nonumber\\
+\mbox{Re}\int_Q (2\lambda^4 {\tau}^3 \varphi^3 a(t, x, \nabla\psi, \nabla\psi)^2
 \vert w\vert^2
 + \vert w\vert^2 \varphi^3 \lambda^3 {\tau}^3 \sum_{{\ell},j=1}^n 
a_{{\ell}j} \partial_{x_j}\psi \partial_{x_{\ell}}a(t, x, \nabla\psi, \nabla\psi)) dxdt\nonumber\\
- \int_\Sigma 2\lambda^3 {\tau}^3  \varphi^3 a(t, x, \nabla\psi, \nabla\psi)
a(t, x, \nabla\psi, \nu)\vert w\vert^2
d\Sigma .
\end{eqnarray}
\end{proposition}
{\bf Proof.} The definition (\ref{sosna}) of the operators $P_1$ and $P_2$ imply
\begin{eqnarray}\label{suka}
\mbox{Re}\, (P_1w,P_2w)_{L^2(Q)}= \mbox{Re}\,  ( i\partial_t w,P_1(t,x,D,\tau)u)_{L^2(Q)} - \mbox{Re}\, (\sum_{{\ell},j=1}^n a_{{\ell}j} \partial^2_{x_{\ell}x_j} w,2 {\tau}\lambda\varphi \sum_{{\ell},j=1}^n
a_{{\ell}j}\partial_{x_{\ell}}\psi\partial_{x_j} w )_{L^2(Q)} \nonumber\\- \mbox{Re}\, (
\lambda^2 {\tau}^2 \varphi^2 a(t, x, \nabla\psi, \nabla\psi) w ,2 {\tau}\lambda\varphi \sum_{{\ell},j=1}^n
a_{{\ell}j}\partial_{x_{\ell}}\psi\partial_{x_j} w )_{L^2(Q)}\nonumber\\
 - \mbox{Re}\, (\sum_{{\ell},j=1}^n a_{{\ell}j} \partial^2_{x_{\ell}x_j} w +
\lambda^2 {\tau}^2 \varphi^2 a(t, x, \nabla\psi, \nabla\psi) w , 2\tau\lambda^2\varphi a(t,x,\nabla\psi,\nabla\psi)u\nonumber\\+2\tau\varphi\lambda \sum_{k,m=1
}^n \partial_{x_k} a_{km}\partial_{x_m}\psi u)_{L^2(Q)}=\sum_{p=0}^3A_p.
\end{eqnarray}
We compute each term $A_j$ separately. We start with term $A_0.$ Integrating by parts  we obtain
$$
(i\partial_t w,P_1(t,x,D,\tau)w)_{L^2(Q)}=( w,P_{1,(0)}(t,x,D,\tau)w)_{L^2(Q)}+( w,iP_1(t,x,D,\tau)\partial_t w )_{L^2(Q)}
$$
$$+
(i w,P_1(t,x,D,\tau)w )_{L^2(\Omega)}\vert_0^T=
$$
$$
( w,P_{1,(0)}(t,x,D,\tau)w)_{L^2(Q)}-(P_1(t,x,D,\tau) w,i\partial_t w )_{L^2(Q)}
$$
$$+
(i w,P_1(t,x,D,\tau)w )_{L^2(\Omega)}\vert_0^T
+ (w, 2{\tau}\lambda i\varphi a(t,x, \nabla\psi, \nu) \partial_tw)_{L^2(\Sigma)}.
$$
Then
\begin{eqnarray}\label{kol}
A_0=\mbox{Re}\,( w,P_{1,(0)}(t,x,D,\tau)w)_{L^2(Q)}+\mbox{Re}\, 
(i w,P_1(t,x,D,\tau)w )_{L^2(\Omega)}\vert_0^T\nonumber\\+\mbox{Re}\, (w, 2{\tau}\lambda i\varphi a(t,x, \nabla\psi, \nu) \partial_tw)_{L^2(\Sigma)}.
\end{eqnarray}

Next, integrating by parts for the second term of right-hand-side of (\ref{suka}) we have
\begin{eqnarray}\label{PD1.56}
A_1  =  \mbox {Re}\,\int_Q - (\sum_{m,j=1}^n a_{mj}
\partial^2_{x_mx_j} w)
(2{\tau}\lambda\varphi \sum_{k,\ell=1}^n a_{k\ell}\partial_{x_k} \psi
\partial_{x_\ell} \bar w) dxdt \nonumber \\
 = \mbox {Re}\, \int_Q
\left(\sum_{m,j=1}^n  \partial_{x_j} a_{mj}
\partial_{x_m} w 2{\tau}\lambda\varphi \sum_{k,\ell=1}^n a_{k\ell}
\partial_{x_k}\psi_{x_k} \partial_{x_\ell}\bar w\right. + 2{\tau}\lambda^2 \varphi \vert a(t, x, \nabla\psi, \nabla w)\vert^2\nonumber \\
 + 2{\tau}\lambda\varphi \sum_{m,j=1}^n a_{{\ell}j}
\partial_{x_m} w \sum_{k,\ell=1}^n
\partial_{x_j} (a_{k\ell}\partial_{x_k} \psi )
\partial_{x_{ \ell}}\bar  w \nonumber \\
 +  2{\tau}\lambda\varphi \sum_{m,j=1}^n a_{mj} \partial_{x_m} w
\sum_{k,\ell=1}^n a_{k\ell}\partial_{x_k} \psi \partial^2_{x_jx_\ell} \bar w\bigg) dxdt
+  \mbox {Re}\, \int_\Sigma 2{\tau}\lambda\varphi 
a(t,x,\nu,\nabla w)a(t,x,\nabla\psi,\nabla\bar w) d\Sigma\nonumber\\
 =  \mbox {Re}\,\int_Q
\left(\sum_{m,j=1}^n \partial_{x_j} a_{mj}
\partial_{x_m} w 2\tau\lambda\varphi \sum_{k,\ell=1}^n a_{k\ell} \psi_{x_k}
\partial_{x_\ell} \bar w + 2{\tau}\lambda^2 \varphi \vert a(t, x, \nabla\psi, \nabla w)\vert^2\right. \nonumber\\
 + 2{\tau}\lambda\varphi \sum_{m,j=1}^n a_{mj} \partial_{x_m} w
\sum_{k, \ell =1}^n \partial_{x_j} (a_{k\ell} \partial_{x_k}\psi)
\partial_{x_\ell}\bar  w
- {\tau}\lambda\varphi
\sum_{k, \ell = 1}^n a_{k\ell} \partial_{x_k}\psi \sum_{m,j=1}^n
\partial_{x_\ell} a_{mj} \partial_{x_m} w
\partial_{x_j} \bar w\nonumber \\
 + \left.{\tau}\lambda \varphi \sum_{k,\ell = 1}^n a_{k\ell}\partial_{x_k} \psi
\partial_{x_\ell} \sum_{m,j=1}^n a_{mj}
\partial_{x_m} w  \partial_{x_j}\bar  w\right) dxdt +
 \mbox {Re}\, \int_\Sigma  2{\tau}\lambda\varphi 
a(t,x,\nu,\nabla w)a(t,x,\nabla\psi,\nabla \bar w) d\Sigma .
\end{eqnarray}

Integrating by parts once again, we obtain from (\ref{PD1.56}):
\begin{eqnarray}\label{D1.57}
A_1 =  \mbox {Re}\, \int_Q \left(\sum_{m,j=1}^n \partial_{x_j} a_{mj}
\partial_{x_m} w2{\tau}\lambda\varphi \sum_{k,\ell=1}^n a_{k\ell} \partial_{x_k}\psi
 \partial_{x_\ell}\bar w\right. + 2{\tau}\lambda^2 \varphi
 \vert a(t, x, \nabla\psi, \nabla w\big)\vert^2+ \nonumber\\
  2{\tau}\lambda \varphi \sum_{m,j=1}^n  a_{mj} \partial_{x_m} w
\sum_{k,\ell=1}^n \partial_{x_j}  (a_{k\ell}\partial_{x_k} \psi)
\partial_{x_\ell} \bar w - {\tau}\lambda\varphi \sum_{k,\ell=1}^n
a_{k\ell} \partial_{x_k}\psi \sum_{m,j=1}^n \partial_{x_\ell} a_{mj}
\partial_{x_m} w \partial_{x_j} \bar w \nonumber\\
 - {\tau}\lambda^2 \varphi a(t, x, \nabla\psi, \nabla\psi) a(t, x, \nabla w, \nabla \bar w)
 - {\tau}\lambda\varphi \sum_{k,\ell=1}^n a_{k\ell}\partial_{x_k} \psi \sum_{m,j=1}^n
\partial_{x_\ell} a_{mj}  \partial_{x_m} w
\partial_{x_j} \bar w\nonumber \\
\left. - a(t, x, \nabla w, \nabla \bar w) {\tau}\lambda\varphi \sum_{k,\ell=1}^n
\partial_{x_\ell} (a_{k\ell}\psi_{x_k})\right) dxdt\nonumber \\
 + \mbox {Re}\,\int_\Sigma \left(
2{\tau}\lambda\varphi 
a(t,x,\nu,\nabla w)a(t,x,\nabla\psi,\nabla\bar w) -
{\tau}\lambda\varphi          a(t, x, \nabla w, \nabla\bar w)
a(t, x, \nu, \nabla\psi)\right) d\Sigma .
\end{eqnarray}

Integrating by parts in the third term of the right-hand-side of (\ref{suka}), we have
\begin{eqnarray}\label{D1.55}
 A_2=- \mbox{Re}\, \int_Q 2\lambda^3 {\tau}^3 w \varphi^3 a(t, x, \nabla\psi, \nabla\psi)
a(t, x, \nabla\psi, \nabla\bar w)   
dxdt \nonumber \\
=  -\int_Q \lambda^3 {\tau}^3 \varphi^3 a(t,x,\nabla\psi,\nabla\psi)
a(t, x, \nabla\psi, \nabla \vert w\vert^2)  dxdt\nonumber \\
=   \int_Q (3\lambda^4 {\tau}^3 \varphi^3 a(t, x, \nabla\psi, \nabla\psi)^2
 \vert w\vert^2
 + \vert w\vert^2 \varphi^3 \lambda^3 {\tau}^3 \sum_{\ell,j=1}^n \partial_{x_\ell}
(a_{\ell j} \partial_{x_j}\psi a(t, x, \nabla\psi, \nabla\psi))) dxdt\nonumber\\
- \int_\Sigma 2\lambda^3 {\tau}^3  \varphi^3 a(t, x, \nabla\psi, \nabla\psi)
a(t, x, \nabla\psi, \nu)\vert w\vert^2
d\Sigma .
\end{eqnarray}
Integrating by parts in the last term of the right-hand-side of (\ref{suka}), we have
\begin{eqnarray}\label{D1.56}
A_3=- \mbox{Re}\, (
\lambda^2 {\tau}^2 \varphi^2 a(t, x, \nabla\psi, \nabla\psi) w , \tau\lambda^2\varphi a(t,x,\nabla\psi,\nabla\psi)u+\tau\varphi\lambda \sum_{k,m=1
}^n \partial_{x_k} a_{km}\partial_{x_m}\psi u)_{L^2(Q)}\nonumber\\
+ \mbox{Re}\, (\sum_{\ell,j=1}^n a_{\ell j} \partial_{x_j} w , \partial_{x_\ell}(\tau\lambda^2\varphi a(t,x,\nabla\psi,\nabla\psi)u+\tau\varphi\lambda \sum_{k,m=1
}^n \partial_{x_k} a_{km}\partial_{x_m}\psi u))_{L^2(Q)}\nonumber\\
- \mbox{Re}\, (\sum_{\ell,j=1}^n a_{\ell j} \partial_{\nu_A} w , \tau\lambda^2\varphi a(t,x,\nabla\psi,\nabla\psi)u+\tau\varphi\lambda \sum_{k,m=1
}^n \partial_{x_k} a_{km}\partial_{x_m}\psi u)_{L^2(\Sigma)}.
\end{eqnarray}
From (\ref{kol}), (\ref{D1.56}), (\ref{D1.57}), (\ref{D1.55}) we obtain (\ref{suka1}). Proof of the proposition is  complete.
$\blacksquare$

\centerline {\bf Proof of Lemma \ref{D1.2}.}

Now we compute function $q_\psi$ explicitly. In order to simplify the notations we set $$p^{(k)}(t,x,\xi)=\partial_{\xi_k}p(t,x,\xi),\quad p_{(k)}(t,x,\xi)=\partial_{x_k}p(t,x,\xi),\quad  p_{(0)}(t,x,\xi)=\partial_{t}p(t,x,\xi).$$

 The short computations imply
$$
\frac{\{p(t,x,\xi-i\tau\nabla_x\psi),p(t,x, \xi+i\tau\nabla_x \psi)\}}{2i\tau}=
$$

$$\frac{1}{2i\tau}\{\sum_{k=0}^n p^{(k)}(t,x,\xi-i\tau\nabla_x\psi)(p_{(k)}(t,x, \xi+i\tau\nabla_x \psi)+\sum_{\ell=1}^np^{(\ell)}(t,x, \xi+i\tau\nabla_x \psi)(i\tau)\partial^2_{x_kx_\ell}\psi)
$$
$$
-\sum_{k=0}^n p^{(k)}(t,x,\xi+i\tau\nabla_x\psi)(p_{(k)}(t,x, \xi-i\tau\nabla_x \psi)-\sum_{\ell=1}^np^{(\ell)}(t,x, \xi-i\tau\nabla_x \psi)(i\tau)\partial^2_{x_kx_\ell}\psi)\}=
$$
$$\frac{1}{2i\tau}( p^{(0)}(t,x,\xi-i\tau\nabla_x\psi)p_{(0)}(t,x, \xi+i\tau\nabla_x \psi)
$$
$$
-p^{(0)}(t,x,\xi+i\tau\nabla_x\psi)p_{(0)}(t,x, \xi-i\tau\nabla_x \psi)+
$$
$$\sum_{k=1}^n p^{(k)}(t,x,\xi-i\tau\nabla_x\psi)(p_{(k)}(t,x, \xi+i\tau\nabla_x \psi)+\sum_{\ell=1}^np^{(\ell)}(t,x, \xi+i\tau\nabla_x \psi)(i\tau)\partial^2_{x_kx_\ell}\psi)
$$
$$
-\sum_{k=1}^n p^{(k)}(t,x,\xi+i\tau\nabla_x\psi)(p_{(k)}(t,x, \xi-i\tau\nabla_x \psi)-\sum_{\ell=1}^np^{(\ell)}(t,x, \xi-i\tau\nabla_x \psi)(i\tau)\partial^2_{x_kx_\ell}\psi))=
$$
$$- 2a_{(0)}(t,x, \xi,\nabla_x \psi)+
$$
$$\frac{1}{2i\tau}(\sum_{k=1}^n p^{(k)}(t,x,\xi-i\tau\nabla_x\psi)(p_{(k)}(t,x, \xi+i\tau\nabla_x \psi)+\sum_{\ell=1}^np^{(\ell)}(t,x, \xi+i\tau\nabla_x \psi)(i\tau)\partial^2_{x_kx_\ell}\psi)
$$
$$
-\sum_{k=1}^n p^{(k)}(t,x,\xi+i\tau\nabla_x\psi)(p_{(k)}(t,x, \xi-i\tau\nabla_x \psi)-\sum_{\ell=1}^np^{(\ell)}(t,x, \xi-i\tau\nabla_x \psi)(i\tau)\partial^2_{x_kx_\ell}\psi))=
$$

$$ -2a_{(0)}(t,x, \xi',\nabla_x \psi)+
\sum_{m,\tilde m,k,\ell=1}^n4 a_{mk}a_{\tilde m \ell}(\xi_m\xi_{\tilde m}+\tau^2\partial_{x_m}\psi\partial_{x_{\tilde m}}\psi) \partial^2_{x_kx_\ell}\psi
$$
$$+
2\sum_{m,k=1}^na_{mk}\partial_{x_m}\psi (a_{(k)}(t,x,\xi',\xi')-\tau^2 a_{(k)}(t,x,\nabla \psi,\nabla\psi))
+4\sum_{m,k=1}^na_{mk}\xi_m a_{(k)}(t,x,\xi',\nabla\psi).$$
Therefore
\begin{eqnarray}\label{somolok1}q_\psi(t,x,\xi,\tau)=
4\sum_{m,\tilde m,k,\ell=1}^na_{mk}a_{\tilde m \ell}(\xi_m\xi_{\tilde m}+\tau^2 \partial_{x_m}\psi\partial_{x_{\tilde m}}\psi) \partial^2_{x_kx_\ell}\psi\nonumber\\-
2\sum_{m,k=1}^na_{mk}\partial_{x_m}\psi (a_{(k)}(t,x,\xi',\xi')-\tau^2 a_{(k)}(t,x,\nabla \psi,\nabla\psi))
+4\sum_{m,k=1}^na_{mk}\xi_m a_{(k)}(t,x,\xi',\nabla\psi)).\end{eqnarray}
Replacing in (\ref{somolok1}) function $\psi$ by $\alpha$ we obtain:
\begin{eqnarray}\label{somolok}q_\alpha (t,x,\xi,\tau)=
4\sum_{m,\tilde m,k,\ell=1}^n \lambda\varphi a_{mk}a_{\tilde m \ell}(\xi_m\xi_{\tilde m}+\tau^2\lambda^2\varphi^2\partial_{x_m}\psi\partial_{x_{\tilde m}}\psi) \partial^2_{x_kx_\ell}\psi\nonumber\\+4\lambda\varphi a(t,x,\nabla\psi,\xi')^2+ 4\tau^2\varphi^2\lambda^4a(t,x,\nabla\psi,\nabla\psi)^2\nonumber\\
-2\lambda\varphi\sum_{m,k=1}^na_{mk}\partial_{x_m}\psi (a_{(k)}(t,x,\xi',\xi')-\tau^2\lambda^2\varphi^2 a_{(k)}(t,x,\nabla \psi,\nabla\psi))
+4\lambda\varphi\sum_{m,k=1}^na_{mk}\xi_m a_{(k)}(t,x,\xi',\nabla\psi)).\nonumber\end{eqnarray}

Since $\nabla \psi\ne 0$ on $\bar Q$ by (\ref{sok5})
we have $a(t,x,\nabla\psi,\nabla \psi)>0$ on $\bar Q.$  Then there exist $\lambda_0$ such that for all $\lambda>\lambda_0$ we have 

\begin{eqnarray}
4\sum_{m,\tilde m,k,\ell=1}^n \lambda\varphi a_{mk}a_{\tilde m \ell}\tau^2\lambda^2\varphi^2\partial_{x_m}\psi\partial_{x_{\tilde m}}\psi \partial^2_{x_kx_\ell}\psi+ 4\tau^2\varphi^2\lambda^4a(t,x,\nabla\psi,\nabla\psi)^2\nonumber\\
+2\lambda\varphi\sum_{m,k=1}^n\tau^2\lambda^2\varphi^2 a_{mk}\partial_{x_m}\psi a_{(k)}(t,x,\nabla \psi,\nabla\psi))>C\tau^2\lambda^2\varphi^3.\nonumber\end{eqnarray}

Then inequality (\ref{pseudo}) holds true if

\begin{eqnarray}\label{somolok1}
4\sum_{m,\tilde m,k,\ell=1}^n \lambda\varphi a_{mk}a_{\tilde m \ell}(\xi_m\xi_{\tilde m}) \partial^2_{x_kx_\ell}\psi\nonumber\\+4\lambda\varphi a(t,x,\nabla\psi,\xi')^2
-2\lambda\varphi\sum_{m,k=1}^na_{mk}\partial_{x_m}\psi (a_{(k)}(t,x,\xi',\xi'))
+4\lambda\varphi\sum_{m,k=1}^na_{mk}\xi_m a_{(k)}(t,x,\xi',\nabla\psi))>0.\nonumber\end{eqnarray}
This proves Proposition \ref{11}.

In particular we have that there exist a strictly positive constant $C$ such that

\begin{equation}\label{pseudo1}
q_\alpha (t,x,\xi,\tau)\ge C(\vert \xi'\vert^2+\tau^2\varphi^2)\quad (t,x,\xi')\in Q\times \Bbb R^{n}, \tau>0.
\end{equation}
 Let us consider the operator
\begin{equation}\label{D1.40}
 \tilde P(t,x,D)u = i\partial_t u - \sum_{{\ell}, j = 1}^n a_{{\ell}j}(t, x)
\partial^2_{x_{\ell} x_j} u.
\end{equation} We set
\begin{equation}\label{D1.41}\tilde{g}(t, x) = g(t, x) - \sum_{{\ell}=1}^n b_{\ell}(t,x) \partial_{x_{\ell}} u - c(t,x)u +
\sum_{{\ell},j = 1}^n \partial_{x_{\ell}} a_{{\ell}j}(t,x) \partial_{x_j} u.
\end{equation}

We denote $w(t,x)=e^{{\tau}\alpha}u(t,x).$
By (\ref{DD1.2}) we have
\begin{equation}\label{D1.42} w(0, \cdot) = 0\quad
\mbox{in}\quad \Omega.
\end{equation}

It follows from (\ref{Ssok1}) and (\ref{D1.40}), (\ref{D1.41}) that
\begin{equation}\label{D1.44}
P(t,x,D,\tau)w  = e^{{\tau}\alpha} \tilde P(t,x,D) e^{-{\tau}\alpha} w
= e^{{\tau}\alpha} \tilde{g} \quad \mbox{in}\quad  Q.
\end{equation}

Operator $P(t,x,D,\tau)$ can be written explicitly as follows
\begin{eqnarray}\label{D1.46}
P(t,x,D,\tau)w = i\partial_t w - \sum_{{\ell},j=1}^n a_{{\ell}j}
\partial^2_{x_{\ell}x_j} w + 2 {\tau}\lambda\varphi \sum_{{\ell},j=1}^n
a_{{\ell}j}\partial_{x_{\ell}}\psi\partial_{x_j} w
  + {\tau}\lambda^2 \varphi
 a(t, x, \nabla\psi, \nabla\psi) w\nonumber\\
 - {\tau}^2\lambda^2\varphi^2  a(t, x, \nabla\psi, \nabla\psi) w
 + {\tau}\lambda \varphi w \sum_{{\ell},j=1}^n a_{{\ell}j} \partial^2_{x_{\ell}x_j}\psi - {i\tau}\partial_t\alpha w.
 \end{eqnarray}
We recall that quadratic form $a(t,x,\xi,\eta)$ was defined in (\ref{pioner}).

We have
\begin{equation}\label{voin}
P_1 w+P_2w=f\quad\mbox{in}\,Q,\quad w\vert_{\Sigma}=0,
\end{equation}
where function $f$ satisfies the estimate
\begin{equation}\label{voin}
\Vert f\Vert_{L^2(Q)}\le C(\Vert ge^{\tau\alpha}\Vert_{L^2(Q)}+\Vert w\Vert_{H^{0,1,\tau}(Q)}).
\end{equation}
Taking the $L^2$ norm of both sides of equation (\ref{voin}) we obtain
\begin{equation}\label{suka22}
\Vert P_1w\Vert^2_{L^2(Q)}+\Vert P_2w\Vert^2_{L^2(Q)}+2\mbox{Re}\,(P_1w,P_2w)_{L^2(Q)}=\Vert f\Vert^2_{L^2(Q)}.
\end{equation}
By proposition \ref{osa}

\begin{eqnarray}\label{Ssuka1}
\mbox{Re} (P_1w,P_2w)_{L^2(Q)}=
\mbox{Re}\, ( w,P_{1,(0)}(t,x,D,\tau)w)_{L^2(Q)}+
\mbox{Re}\, (i w(T,\cdot),P_1(T,x,D,\tau)w (T,\cdot))_{L^2(\Omega)} \nonumber\\
   +\mbox{Re}\int_Q \left(\right.  2{\tau}\lambda^2 \varphi
 \vert a(t, x, \nabla\psi, \nabla w\big)\vert^2 \nonumber\\
\left.
  +2{\tau}\lambda \varphi \sum_{m,j=1}^n a_{mj} \partial_{x_m} w
\sum_{k,\ell=1}^n \partial_{x_j}  (a_{k\ell}\partial_{x_k} \psi)
\partial_{x_\ell} \bar w - {\tau}\lambda\varphi \sum_{k,\ell=1}^n
a_{k\ell} \partial_{x_k}\psi \sum_{m,j=1}^n \partial_{x_\ell} a_{mj}
\partial_{x_m} w \partial_{x_j} \bar w \right) dxdt\nonumber \\
 + \int_\Sigma 
{\tau}\lambda\varphi 
a(t,x,\nu,\nu)a(t,x,\nabla\psi,\nu)\vert\partial_\nu w\vert^2 
 d\Sigma \nonumber\\
+\int_Q (2\lambda^4 {\tau}^3 \varphi^3 a(t, x, \nabla\psi, \nabla\psi)^2
 \vert w\vert^2
 + \vert w\vert^2 \varphi^3 \lambda^3 {\tau}^3 \sum_{\ell,j=1}^n 
a_{{\ell}j} \partial_{x_j}\psi \partial_{x_{\ell}}a(t, x, \nabla\psi, \nabla\psi)) dxdt.
\end{eqnarray}
The following is true

\begin{proposition}\label{supostat} {\it 
Let $z=(t_0,y) $ be an arbitrary point from $\bar Q.$ There exist constants $\delta=\delta(z), K(z)>0, C(z)>0, \tilde C(z)>0$ such that  if $\mbox{supp}\subset B(z,\delta)$ then
\begin{eqnarray} \label{Z}\frac{K}{\tau}\Vert P_1 w\Vert^2_{L^2(B(z,\delta))}
+\mbox{Re}\int_{B(z,\delta)} \left(\right.  2{\tau}\lambda^2 \varphi
 \vert a(t, x, \nabla\psi, \nabla w\big)\vert^2 \nonumber\\
  +2{\tau}\lambda \varphi \sum_{m,j=1}^n a_{mj} \partial_{x_m} w
\sum_{k,\ell=1}^n \partial_{x_j}  (a_{k\ell}\partial_{x_k} \psi)
\partial_{x_\ell} \bar w - {\tau}\lambda\varphi \sum_{k,\ell=1}^n
a_{k\ell} \partial_{x_k}\psi \sum_{m,j=1}^n \partial_{x_\ell} a_{mj}
\partial_{x_m} w \partial_{x_j} \bar w ) dxdt\nonumber\\
+\int_{B(z,\delta)} (2\lambda^4 {\tau}^3 \varphi^3 a(t, x, \nabla\psi, \nabla\psi)^2
 \vert w\vert^2
 + \vert w\vert^2 \varphi^3 \lambda^3 {\tau}^3 \sum_{{\ell},j=1}^n 
a_{{\ell}j} \partial_{x_j}\psi \partial_{x_{\ell}}a(t, x, \nabla\psi, \nabla\psi)) dxdt\ge\nonumber\\ C(z)\int_{B(z,\delta)} (\tau\varphi \vert \nabla w\vert^2+\tau^3\varphi^3\vert w\vert^2)dxdt-\tilde C(z)\int_{B(z,\delta)} (\vert \nabla w\vert^2+\tau^2\varphi^2\vert w\vert^2)dxdt\quad \forall \tau\ge 0.
\end{eqnarray}
}
\end{proposition}
{\bf Proof.} For any $\epsilon>0$ there exist $\delta(z,\epsilon)$ such that 
\begin{eqnarray}\label{Picho1}\vert \mbox{Re}\,
\int_{B(z,\delta)} \left(\right.  2{\tau}\lambda^2 \varphi
 \vert a(t, x, \nabla\psi, \nabla w\big)\vert^2 \nonumber\\
  +2{\tau}\lambda \varphi \sum_{m,j=1}^n a_{mj} \partial_{x_m} w
\sum_{k,\ell=1}^n \partial_{x_j}  (a_{k\ell}\partial_{x_k} \psi)
\partial_{x_\ell} \bar w - {\tau}\lambda\varphi \sum_{k,m=1}^n
a_{km} \partial_{x_k}\psi \sum_{m,j=1}^n \partial_{x_\ell} a_{mj}
\partial_{x_m} w \partial_{x_j} \bar w ) dxdt\nonumber\\
+\mbox{Re}\,\int_{B(z,\delta)} (2\lambda^4 {\tau}^3 \varphi^3 a(z, \nabla\psi, \nabla\psi)^2
 \vert w\vert^2
 + \vert w\vert^2 \varphi^3\lambda^3 {\tau}^3 \sum_{\ell,j=1}^n 
a_{\ell j}(z) \partial_{x_j}\psi \partial_{x_\ell}a(z, \nabla\psi, \nabla\psi)) dxdt\nonumber\\-\mbox{Re}\,\int_{B(z,\delta)} \left(\right.  2{\tau}\lambda^2 \varphi(t,z)
 \vert a(z, \nabla\psi(y), \nabla w\big)\vert^2 \nonumber\\
  +2{\tau}\lambda \varphi(t,y) \sum_{\ell,j=1}^n a_{\ell j}(z) \partial_{x_\ell} w
\sum_{k,m=1}^n \partial_{x_j}  (a_{km}(z)\partial_{x_k} \psi(y))
\partial_{x_m} \bar w\nonumber\\ - {\tau}\lambda\varphi(t,y) \sum_{k,\ell=1}^n
a_{k\ell}(z) \partial_{x_k}\psi(y) \sum_{m,j=1}^n \partial_{x_\ell} a_{mj}(z)
\partial_{x_m} w \partial_{x_j} \bar w ) dxdt\nonumber\\
+\mbox{Re}\,\int_{B(z,\delta)} (2\lambda^4 {\tau}^3 \varphi^3(t,z) a(z, \nabla\psi(y), \nabla\psi(y))^2
 \vert w\vert^2\nonumber\\
 + \vert w\vert^2 \varphi^3(t,y) \lambda^3 {\tau}^3 \sum_{\ell,j=1}^n 
a_{\ell j}(z) \partial_{x_j}\psi(y) \partial_{x_i}a(z, \nabla\psi(y), \nabla\psi(y))) dxdt\vert\nonumber\\\le \epsilon \int_{B(z,\delta)} (\tau\varphi \vert \nabla w\vert^2+\tau^3\varphi^3\vert w\vert^2)dxdt
\end{eqnarray}

and
\begin{equation}\label{Picho2}
\vert\frac{K}{\tau}\Vert P_1 w\Vert^2_{L^2(B(z,\delta))}-\Vert \frac{K}{\tau\varphi(t,y)} \tilde P_1 w\Vert^2_{L^2(B(z,\delta))}\vert \le \epsilon \int_{B(z,\delta)} (\tau\varphi \vert \nabla w\vert^2+\tau^3\varphi^3\vert w\vert^2)dxdt,
\end{equation}
where
$$
\tilde P_1(t,x,D,\tau)u=2\tau\varphi(t,y)\lambda a(z,\nabla \psi(y),\nabla u)+\tau\lambda^2\varphi(t,y) a(z,\nabla\psi(y),\nabla\psi(y))u+\tau\varphi(t,y)\lambda \sum_{k,j=1
}^n \partial_{x_k} a_{kj}(z)\partial_{x_j}\psi(y) u .
$$
Next for each $t\in (0,T)$ we extend the function $u$ by zero outside of domain $\Omega.$ We set 
\begin{eqnarray}\label{sifilis}
\mathcal I(w)=\int_{[0,T]\times \Bbb R^n} (2\lambda^4 {\tau}^3 \varphi^3 a(z, \nabla\psi, \nabla\psi)^2
 \vert w\vert^2
 + \vert w\vert^2 \varphi^3\lambda^3 {\tau}^3 \sum_{\ell ,j=1}^n 
a_{\ell j}(z) \partial_{x_j}\psi \partial_{x_\ell}a(z, \nabla\psi, \nabla\psi)) dxdt\nonumber\\-\int_{[0,T]\times \Bbb R^n} \left(\right.  2{\tau}\lambda^2 \varphi(t,z)
 \vert a(z, \nabla\psi(y), \nabla w\big)\vert^2
  +2{\tau}\lambda \varphi(t,y) \sum_{\ell,j=1}^n a_{\ell j}(z) \partial_{x_\ell} w
\sum_{k,\ell=1}^n \partial_{x_j}  (a_{k\ell}(z)\partial_{x_k} \psi(y))
\partial_{x_m} \bar w  \nonumber\\- {\tau}\lambda\varphi(t,y) \sum_{k,\ell=1}^n
a_{km}(z) \partial_{x_k}\psi(y) \sum_{\ell,j=1}^n \partial_{x_m} a_{\ell j}(z)
\partial_{x_\ell} w \partial_{x_j} \bar w ) dxdt\nonumber\\
+\int_{[0,T]\times \Bbb R^n}(2\lambda^4 {\tau}^3 \varphi^3(t,z) a(z, \nabla\psi(y), \nabla\psi(y))^2
 \vert w\vert^2\nonumber\\
 + \vert w\vert^2 \varphi^3(t,y) \lambda^3 {\tau}^3 \sum_{\ell,j=1}^n 
a_{\ell j}(z) \partial_{x_j}\psi(y) \partial_{x_\ell}a(z, \nabla\psi(y), \nabla\psi(y))) dxdt+\frac{K}{\tau}\Vert \tilde P_1 w\Vert^2_{L^2([0,T]\times \Bbb R^n))}.
\end{eqnarray}
Taking the Fourier transform respect to the variable $x$ in (\ref{sifilis}) we obtain
$$
\mathcal I (w)
=\int_{[0,T]\times \Bbb R^n}\tau \varphi(t,y)( q_\alpha (t,\xi',\tau\varphi(t,y))+4\lambda^2 K a(z,\nabla \psi(y),\xi'))\vert Fw\vert^2d\xi'dt.
$$
Taking constant $K$ sufficiently large, by (\ref{pseudo1}), we obtain that  there exist some positive $\tilde C$ such that
$$
q_\alpha (t,\xi',\tau\varphi(t,y))+4\lambda^2 K a(z,\nabla \psi(y),\xi')\ge \tilde C(\vert \xi'\vert^2+\tau^2/t^2)\quad \forall( t,\xi',\tau)\in [0,T]\times \Bbb R^n\setminus\{0\}\times \ R^1_+.
$$
Hence
$$
\mathcal I(w)\ge \tilde C\int_{Q}(\tau\varphi\vert \nabla w\vert^2+\tau^3\varphi^3\vert w\vert^2)dxdt.
$$
Taking in (\ref{Picho1}), (\ref{Picho2}) parameter $\epsilon$ sufficiently small we obtain (\ref{Z}). Proof of the proposition is complete.$\blacksquare $

Next we drop the assumption $\mbox{supp}\, w\subset B(z,\delta)$ in the proposition \ref{supostat}.
From the covering of the set $\bar Q$ by balls $B(z,\delta(z))$ we will take the finite subcovering $B(z_k,\delta_k)$ with $k\in\{1,\dots,N\}.$ Let $e^2_k$ be a partition of unity subjected to this covering:
$$
e_k\in C^\infty_0(B(z_k,\delta_k)), \quad e_k\ge 0, \quad \sum_{k=1}^Ne_k^2=1\quad \mbox{on}\quad \Omega.
$$ We set $w_k=e^2_kw, \tilde w_k=e_k w$. Then $w=\sum_{k=1}^Nw_k.$ The functional $\mathcal I$ can be written in the form
$$
\mathcal I(w)=\mbox{Re}\, \sum_{\ell=1}^M (K_{1,\ell}(x,D,\tau) w,K_{2,\ell} (x,D,\tau)w)_{L^2(Q)}.
$$
Then 
$$
\mathcal I(w)=\mbox{Re}\, \sum_{\ell=1}^M(K_{1,\ell}(x,D,\tau) w,K_{2,\ell}(x,D,\tau)w)_{L^2(Q)}=\sum_{k=1}^N \mbox{Re}\, (K_{1,\ell}(x,D,\tau) w,K_{2,\ell}(x,D,\tau)w_k)_{L^2(Q)}=
$$
$$\sum_{\ell=1}^M\sum_{k=1}^N \mbox{Re}\, (K_{1,\ell}(x,D,\tau) w,K_{2,\ell}(x,D,\tau)(e_k\tilde w_k))_{L^2(Q)}=\sum_{\ell=1}^M\sum_{k=1}^N (\mbox{Re}\, (K_{1,\ell}(x,D,\tau) w,e_kK_{2,\ell}(x,D,\tau)\tilde w_k)_{L^2(Q)}
$$
$$+\mbox{Re}\, (K_{1,\ell}(x,D,\tau) w,[e_k,K_{2,\ell}]\tilde w_k)_{L^2(Q)})=
$$
$$\sum_{\ell=1}^M\sum_{k=1}^N (\mbox{Re}\, (K_{1,\ell}(x,D,\tau) \tilde w_k,K_{2,\ell}(x,D,\tau)\tilde w_k)_{L^2(Q)}+\mbox{Re}\, ([e_1,K_{1,\ell}] w,K_{2,\ell}(x,D,\tau)\tilde w_k)_{L^2(Q)}
$$
$$+\mbox{Re}\, (K_{1,\ell}(x,D,\tau) w,[e_k,K_{2,\ell}]\tilde w_k)_{L^2(Q)}).
$$
Applying the proposition \ref{supostat} we obtain
$$
\mathcal I(w)\ge\sum_{k=1}^N\int_{B(z_k,\delta)} C(z_k)(\tau\varphi \vert \nabla \tilde w_k\vert^2+\tau^3\varphi^3\vert \tilde w_k\vert^2)dxdt-\tilde C(z_k)\int_{B(z_k,\delta)} (\vert \nabla \tilde w_k\vert^2+\tau^2\varphi^2\vert \tilde w_k\vert^2)dxdt
$$
$$+\sum_{\ell=1}^M\mbox{Re}\, ([e_1,K_{1,\ell}] w,K_{2,\ell}(x,D,\tau)\tilde w_k)_{L^2(Q)}+\mbox{Re}\, (K_{1,\ell}(x,D,\tau) w,[e_k,K_{2,\ell}]\tilde w_k)_{L^2(Q)}.
$$ 
 Let $C=\mbox{min}_{k\in\{1,\dots,N\}}\{C(z_k)\}.$ From the above inequality we have
$$
\mathcal I(w)\ge \sum_{k=1}^N\int_{B(z,\delta)} C(\tau\varphi ((\nabla w,\nabla w_k)+(w\nabla e_k,\nabla\overline{\tilde w_k})+(w,\overline {\tilde w_k}\nabla e_k))+\tau^3\varphi^3(w, w_k))dxdt
$$
$$-\tilde C(z_k)\int_{B(z,\delta)} (\vert \nabla \tilde w_k\vert^2+\tau^2\varphi^2\vert \tilde w_k\vert^2)dxdt
$$
$$+\mbox{Re}\, ([e_1,K_{1,\ell}] w,K_{2,\ell}(x,D,\tau)\tilde w_k)_{L^2(Q)}+\mbox{Re}\, (K_{1,\ell}(x,D,\tau) w,[e_k,K_{2,\ell}]\tilde w_k)_{L^2(Q)}
$$
$$
+C\int_{Q} (\tau\varphi \vert\nabla w\vert^2+\tau^3\varphi^3\vert w\vert^2)dxdt
$$

$$+\mbox{Re}\, ([e_1,K_{1,\ell}] w,K_{2,\ell}(x,D,\tau)\tilde w_k)_{L^2(Q)}+\mbox{Re}\, (K_{1,\ell}(x,D,\tau) w,[e_k,K_{2,\ell}]\tilde w_k)_{L^2(Q)}.
$$
By  the Cauchy inequality we obtain 
$$
\mathcal I(w)=\frac C2\int_{Q} (\tau\varphi \vert\nabla w\vert^2+\tau^3\varphi^3\vert w\vert^2)dxdt-\tilde C\int_{Q} (\vert \nabla w\vert^2+\tau^2\varphi^2\vert w\vert^2)dxdt.
$$

From proposition \ref{supostat}, equality (\ref{suka1}) and (\ref{suka22}) we have
\begin{eqnarray} 
\int_Q(\tau \varphi\vert\nabla w\vert^2+\tau^3\varphi^3\vert w\vert^2)dxdt+\sum_{k=1}^2\Vert P_kw\Vert^2_{L^2(Q)}\nonumber\\\le C(\Vert f\Vert^2_{L^2(Q)}+\vert  \int_{\Sigma_0} 
{\tau}\lambda\varphi 
a(t,x,\nu,\nu)a(t,x,\nabla\psi,\nu)\vert\partial_\nu w\vert^2 
 d\Sigma\vert +\vert \mbox{Re}\, (i w(T,\cdot),P_1(T,x,D,\tau)w (T,\cdot))_{L^2(\Omega)}\vert\nonumber\\
-C\int_{Q} (\vert \nabla w\vert^2+\tau^2\varphi^2\vert w\vert^2)dxdt\quad \forall \tau\ge 0.\end{eqnarray}
Using the estimate (\ref{voin}) we obtain

\begin{eqnarray} 
\int_Q(\tau \varphi\vert\nabla w\vert^2+\tau^3\varphi^3\vert w\vert^2)dxdt+\sum_{k=1}^2\Vert P_kw\Vert^2_{L^2(Q)}\nonumber\\\le C(\Vert ge^{\tau\alpha}\Vert^2_{L^2(Q)}+\vert  \int_{\Sigma_0} 
{\tau}\lambda\varphi 
a(t,x,\nu,\nu)a(t,x,\nabla\psi,\nu)\vert\partial_\nu w\vert^2 
 d\Sigma\vert +\vert \mbox{Re}\, (i w(T,\cdot),P_1(T,x,D,\tau)w (T,\cdot))_{L^2(\Omega)}\vert)\end{eqnarray}
for all $\tau\ge\tau_1.$
Proof of lemma is complete. $\blacksquare$

\end{document}